\documentclass[final,leqno,onefignum,onetabnum]{siamltex}

\usepackage{algorithmic,varwidth}
\usepackage{setspace}
\usepackage{amsmath, amssymb, lscape}
\usepackage{enumerate}
\usepackage{amsfonts}
\usepackage{graphicx}
\usepackage{rotating,graphics,psfrag}
\usepackage[dvips]{epsfig}
\usepackage{verbatim}
\usepackage{float}
\usepackage{appendix}
\usepackage{amsbsy}
\usepackage{multirow}
\usepackage{url}
\usepackage{arydshln}
\usepackage{setspace}
\newcommand{\ii}{\i}

\title{Parallel Solution of the Linear Elasticity problem with
applications in Topology Optimisation}

\author{J. Turner\footnotemark[3]
\and M. Ko\v{c}vara\footnotemark[2]
\and D. Loghin\footnotemark[3]}

\newcommand*{\defeq}{\mathrel{\vcenter{\baselineskip0.5ex \lineskiplimit0pt\hbox{\small.}\hbox{\small.}}}=}
\newcommand*{\rdefeq}{=\mathrel{\vcenter{\baselineskip0.5ex \lineskiplimit0pt
                     \hbox{\small.}\hbox{\small.}}}%
                     }

\begin{document}
\maketitle
\renewcommand{\thefootnote}{\fnsymbol{footnote}}

\footnotetext[2]{School of
    Mathematics, University of Birmingham, Edgbaston, Birmingham B15 2TT, UK,
    and Institute of Information Theory
and Automation, Academy of Sciences of the Czech Republic, Pod
vod\'arenskou v\v{e}\v{z}\'{\ii}~4, 18208 Praha 8, Czech Republic. The
work of this author has been partly supported by the EU FP7 project
AMAZE and by grant A100750802 of the Czech Academy of Sciences}
\footnotetext[3]{School of
    Mathematics, University of Birmingham, Edgbaston, Birmingham B15 2TT, UK}

\renewcommand{\thefootnote}{\arabic{footnote}}

\begin{abstract}
In this paper, we aim to solve the system of equations governing linear elasticity in parallel using domain decomposition. Through a non-overlapping decomposition of the domain, our approach aims to target the resulting interface problem, allowing for the parallel computation of solutions in an efficient manner. As a major application of our work, we apply our results to the field of topology optimisation, where typical solvers require repeated solutions of linear elasticity problems resulting from the use of a Picard approach.
\end{abstract}

\begin{keywords}Linear elasticity, topology optimisation, domain decomposition, preconditioning, Krylov methods \end{keywords}


\pagestyle{myheadings} \thispagestyle{plain} \markboth{Parallel Solution of the Linear Elasticity problem with
applications in Topology Optimisation
}{Parallel Solution of the Linear Elasticity problem with
applications in Topology Optimisation
}

\section{Introduction}

Consider a solid elastic body occupying an open and connected domain $\Omega \subset \mathbb{R}^d$ with Lipschitz boundary $\partial\Omega = \partial\Omega_D \cup \partial\Omega_N$, where clamping and traction are imposed on $\partial\Omega_D$ and $\partial\Omega_N$ respectively. Under the application of both body forces $f \colon \Omega \rightarrow \mathbb{R}^d$ and boundary tractions $g \colon \partial\Omega_N \rightarrow \mathbb{R}^d$ the material is subject to deformation so that a given reference point $\mathbf{x}$ of the initial undeformed material is translated to the vector $\mathbf{x}' = \mathbf{x} + u(\mathbf{x})$ of the deformed material, with $u$ denoting the displacement. Through the assumption of linearly elastic material behaviour, the governing equations for $u$ correspond to the following mixed boundary value problem
\begin{subequations}
\label{ElastProblem}
\begin{align}
\hspace{0.6in} \mathcal{L} u \defeq \nabla \cdot \pmb{\sigma} (u)  &= f, \quad  &\textrm{in} \, \, \, &\Omega,  \\ 
\hspace{0.6in}\pmb{\sigma} (u) &= E \colon \pmb{\epsilon} (u), \quad \hspace{0.1in} &\textrm{in}\, \, &\Omega, \label{ElastProblemb} \\ 
\hspace{0.6in}u &= 0, \quad  &\textrm{on}\, \, &\partial\Omega_D, \\
\hspace{0.6in}\pmb{\sigma} (u) \cdot \hat{\mathbf{n}} &= g, \quad  &\textrm{on}\, \, &\partial\Omega_N, \end{align}
\end{subequations}
In the above, the strain caused as a result of the displacements $u$ is characterised by the symmetric linearised strain tensor 
\begin{align*}
\pmb{\epsilon} (u) = \{ \epsilon_{i j} (u) \}_{i,j = 1}^{d} , \hspace{0.5in} \epsilon_{i j} (u) = \frac{1}{2}\left(\frac{\partial u_i}{\partial x_j} + \frac{\partial u_j}{\partial x_i}\right).  
\end{align*}
Additionally, $\hat{\mathbf{n}}$ corresponds to the unit outward pointing normal vector on $\partial\Omega$ and $E \defeq E(\mathbf{x})$ denotes the fourth order elasticity tensor, describing the elastic stiffness of $\Omega$ as a result of the load placed upon it.

We will consider the case where our body consists of one or more isotropic materials (i.e: rotational and directional independence). Equation (\ref{ElastProblemb}) describing Hooke's Law can be written as
\[
\pmb{\sigma} (u) = 2\mu \, \pmb{\epsilon} (u) + \hat{\lambda} \, \text{tr}\left(\pmb{\epsilon} (u) \right) \hat{I},
\]
where $\hat{I}$ represents the identity matrix of appropriate size, $\text{tr}(A)$ denotes the trace of a matrix $A$ and both $\mu$ and $\hat{\lambda}$ correspond to Lam\'e constants defined in the usual manner
\[
\mu \defeq \frac{\bar{E}}{2(1+\nu)}, \hspace{0.5in} \hat{\lambda}\defeq \frac{\nu \bar{E}}{(1+\nu)(1-2\nu)},
\]
with $\bar{E}>0$ corresponding to Young's modulus and $-1 < \nu < 1/2$ the Poisson ratio. We now look to apply domain decomposition to the problem (\ref{ElastProblem}). To do this, we divide our domain $\Omega$ into $N$ nonoverlapping subdomains $\Omega_i$ with local boundaries $\partial\Omega_i$ and outer unit normals $\hat{\mathbf{n}}_i$. We denote by $\Gamma$ the resulting skeletal interface $\Gamma = \cup^{N}_{i=1} \Gamma_i$ where $\Gamma_i \defeq \partial\Omega_i \backslash \partial\Omega$ and by $I \defeq \bar{\Omega}\backslash \Gamma$ the set of interior nodes, with $u_{|_{\Omega_i}}\rdefeq u_i$. Assuming that the restriction of $u_i$ to components of the skeletal interface $u_{i_{|_{\Gamma_i}}}$ is known, problem (\ref{ElastProblem}) is equivalent to the following set of subproblems
\vspace{-0.1in}
\begin{spacing}{1.1}
\begin{align}
\begin{cases}
\begin{aligned}
\label{DDElastProblem}
\mathcal{L} u_i  &= f_i, \quad  &\textrm{in} \, \, \, &\Omega_i, \hspace{0.3in} \\ 
u_i &= 0, \quad \hspace{0.1in} &\textrm{on}\, \, &\partial\Omega_D \cap \partial\Omega_i, \hspace{0.3in} \\ 
\pmb{\sigma} \left( u_i \right) \cdot \mathbf{n}_i &= g_i, \quad  &\textrm{on}\, \, &\partial\Omega_N \cap \partial\Omega_i, \hspace{0.3in} \\
u_i &= \lambda_i, \quad  &\textrm{on}\, \, &\Gamma_i, \hspace{0.3in} \end{aligned}
\end{cases}
\end{align}
\end{spacing}
\vspace{-0.01in}
\noindent where $i=1,\dots,N$. By writing $u_i = u^{(1)}_i + u^{(2)}_i$, we look to describe an appropriate interface operator that will allow problems to be decoupled and thus solved strictly on subdomains in parallel. Through the definition of matrix extension operators $H_i$ that map interface data to relevant subdomains via $u^{(2)}_i = H_i \lambda$, the system (\ref{DDElastProblem}) can be decoupled into the following $2N+1$ subproblems
\begin{subequations}
\label{DDElasProb}
\begin{align}
&\begin{cases}
\label{SplitSigma1}
\begin{aligned}
\mathcal{L} u^{(1)}_i &= f_i, \quad  &\textrm{in} \, \, \, &\Omega_i, \hspace{0.3in} \\ 
u^{(1)}_i &= 0, \quad \hspace{0.1in} &\textrm{on}\, \,& \partial\Omega_D \cap \partial\Omega_i, \hspace{0.3in}  \\
\pmb{\sigma} ( u^{(1)}_i ) \cdot \mathbf{n}_i &= g_i, \quad  &\textrm{on}\, \,& \partial\Omega_N \cap \partial\Omega_i, \hspace{0.3in} \\
u^{(1)}_i &= 0, \quad  &\textrm{on}\, \, &\Gamma_i. \hspace{0.3in} 
\end{aligned}
\end{cases}
\\
&\begin{cases}
\label{2ndprobsplit}
\begin{aligned}
\displaystyle\sum_{i=1}^N  \pmb{\sigma}(H_i \lambda) \cdot \mathbf{n}_i &= - \displaystyle\sum_{i=1}^N  \pmb{\sigma}(u_{i}^{(1)}) \cdot \mathbf{n}_i,  \,\, & &\textrm{on} \,\, \Gamma.  
\end{aligned}
\end{cases}
\\
&\begin{cases}
\label{SplitSigma2}
\begin{aligned}
\mathcal{L} u^{(2)}_i &= 0, \quad  &\textrm{in} \, \, \, &\Omega_i, \hspace{0.3in} \\ 
u^{(2)}_i &= 0, \quad \hspace{0.1in} &\textrm{on}\, \,& \partial\Omega_D \cap \partial\Omega_i, \hspace{0.3in}  \\
\pmb{\sigma} ( u^{(2)}_i ) \cdot \mathbf{n}_i &= 0, \quad  &\textrm{on}\, \,& \partial\Omega_N \cap \partial\Omega_i, \hspace{0.3in} \\
u^{(2)}_i &= \lambda_i, \quad  &\textrm{on}\, \, &\Gamma_i. \hspace{0.3in} 
\end{aligned}
\end{cases}
\end{align}
\end{subequations}
The associated weak form to problem (\ref{2ndprobsplit}) is referred to as the Steklov-Poincar\'e equation, where the so-called Steklov-Poincar\'e pseudo-differential operator $S \colon \Lambda_{\theta} \rightarrow \Lambda_{\theta}'$ defined in the following manner \cite{Turner14}
\[
s(\lambda, \eta) = \langle S \lambda, \eta \rangle \defeq \displaystyle\sum_{i=1}^N \left[ \displaystyle\int_{\Gamma_i} \left( \pmb{\sigma} \left( H_i \lambda_i \right) \cdot \mathbf{n}_i \right) \eta_i \,\text{d}s \right] \rdefeq \displaystyle\sum_{i=1}^N \langle S_i \lambda_i, \eta_i \rangle,
\]
where $\lambda, \eta \in \Lambda_{\theta}$ and $\lambda_{|_{\Gamma_i}} \rdefeq \lambda_i$, $\eta_{|_{\Gamma_i}} \rdefeq \eta_i$. The space $\Lambda_{\theta}$ is chosen to be a suitable fractional Sobolev space of index $\theta$ based on the boundary conditions of the problem, dependent on the intersection of $\Gamma$ with $\partial\Omega$ \cite{Quarteroni99, Toselli05}.
  
\section{Matrix Formulation}

Through a finite element discretisation of the weak formulation, it can be shown that the discrete formulation to the original problem (\ref{ElastProblem}) requires the solution of a matrix-vector system \cite{Turner14}. By distributing nodes based on their location within the domain, we can view this system as follows
\begin{equation}
\label{NumberFour}
K\mathbf{u} = \begin{pmatrix} K_{II} & K_{I \Gamma} \\ K_{\Gamma I} & K_{\Gamma \Gamma} \end{pmatrix} \begin{pmatrix} \mathbf{u}_I \\ \mathbf{u}_{\Gamma} \end{pmatrix} = \begin{pmatrix} \mathbf{f}_I \\ \mathbf{f}_{\Gamma} \end{pmatrix} = \mathbf{f},
\end{equation}
where
\begin{equation}
\label{NumberFive}
K_{II} \defeq \bigoplus_{i=1}^N K_{I_i I_i}, \hspace{0.5in}  \mathbf{u} = \left( \mathbf{u}_I, \mathbf{u}_{\Gamma} \right)^T \in \mathbb{R}^{n=n_I+n_{\Gamma}}.
\end{equation}
In comparison, the corresponding matrix formulations for each of the discrete weak formulations to the problems presented in (\ref{DDElastProblem}) can be written down as
\vspace{-0.1in}
\begin{spacing}{1.1}
\begin{align}
\begin{cases}
\label{3MatSplit}
\begin{aligned}
K_{II} \mathbf{u}^{(1)}_I &= \mathbf{f}_{I},\\
S\mathbf{u}_{\Gamma} &= \mathbf{f}_{\Gamma} - K_{\Gamma I} \mathbf{u}^{(1)}_I,\\
K_{II} \mathbf{u}^{(2)}_I &= - K_{I\Gamma} \mathbf{u}_{\Gamma},
\end{aligned}
\end{cases}
\end{align}
\end{spacing}
\vspace{-0.01in}
\noindent with global solution $\mathbf{u} = \left( \mathbf{u}_I^{(1)},  \mathbf{0} \right)^T + \left( \mathbf{u}_I^{(2)},  \mathbf{u}_{\Gamma} \right)^T$. In the above, the matrix $S$ corresponds to the Schur complement, and so the discretisation of the decoupled problem (\ref{DDElastProblem}) can be viewed as a Schur complement approach to the discretisation of the global problem (\ref{ElastProblem}). Using (\ref{NumberFour}) and (\ref{NumberFive}), we are able to view (\ref{3MatSplit}) in terms of $2N+1$ subproblems in the following way
\begin{align}
\begin{cases}
\label{3MatSplit2}
\begin{aligned}
K_{I_i I_i} \mathbf{u}^{(1)}_{I_i} &= \mathbf{f}_{I_i}, \hspace{1.13in} i =1, \dots, N,\\
\hspace{0.17in} S\mathbf{u}_{\Gamma} &= \mathbf{f}_{\Gamma} - \displaystyle\sum_{i=1}^N K_{\Gamma I_i} \mathbf{u}^{(1)}_{I_i},\\
K_{I_i I_i} \mathbf{u}^{(2)}_{I_i} &= - K_{I_i \Gamma} \mathbf{u}_{\Gamma}, \hspace{0.72in} i =1, \dots, N.
\end{aligned}
\end{cases}
\end{align}
We therefore look to solve (\ref{NumberFour}) by exploiting the potential for parallelisation present in (\ref{3MatSplit2}).

\section{Preconditioning}
\label{SectionThree}
The systems we expect to solve will typically be both sparse and large scale, due to the expected fineness of the finite element discretisation required in modern design processes, allowing for the computation of resolute solutions. This is of particular importance for domains containing sharp jumps, occurring for instance due to predefined fixed or void regions. Therefore, it is appropriate to consider iterative solution techniques when solving systems of the form (\ref{NumberFour}). For our problem, we will consider GMRES \cite{Saad86} for reasons to be described below. 

In order to avoid the direct construction and application of the Schur complement matrix, and to improve the spectral properties of the system matrix, we seek an appropriate preconditioner for the system (\ref{NumberFour}). Through the following choice of $P$, we see that
\[
P=\begin{pmatrix} K_{II} & K_{I\Gamma}\\0 & S \end{pmatrix}, \hspace{0.5in} K P^{-1}=\begin{pmatrix} \hat{I}_{II} & 0\\K_{\Gamma I} K_{II}^{-1} & \hat{I}_{\Gamma\Gamma} \end{pmatrix}
\]
The minimum polynomial of $K P^{-1}$ is $(\lambda - 1)^{d}$,  suggesting that iterative solution methods such as GMRES will converge in at most $d$ iterations \cite{Ipsen02}. Based on this, we propose to precondition from the right with an approximation $\widetilde{P}$ of $P$ as follows
\[
K \widetilde{P}^{-1} \tilde{\mathbf{u}} = \mathbf{f}, \hspace{0.5in} \tilde{\mathbf{u}} = \widetilde{P} \mathbf{u},
\]
where
\[
\widetilde{P}=\begin{pmatrix} K_{II} & K_{I\Gamma}\\0 & \widetilde{S} \end{pmatrix}, \hspace{0.5in} \widetilde{P}^{-1} = \begin{pmatrix} K_{II}^{-1} & 0\\0 & \hat{I}_{\Gamma\Gamma} \end{pmatrix}\begin{pmatrix} \hat{I}_{II} & -K_{I\Gamma}\\0 & \hat{I}_{\Gamma\Gamma} \end{pmatrix}\begin{pmatrix} \hat{I}_{II} & 0\\0 & \widetilde{S}^{-1} \end{pmatrix}
\]
with $\widetilde{S}$ representing an approximation to the discrete Steklov-Poincar\'e operator. We therefore seek a representation of $\widetilde{S}$ ̃that is not only practical to invert but can also be seen to provide an appropriate preconditioning strategy for the resulting interface problem.

The form of $\widetilde{S}$ chosen is based on work in \cite{Arioli09}, where discrete norm representations for projections of the interpolation spaces $\Lambda_{\theta}$ onto suitable finite dimensional subspaces are described and analysed. Discrete norms of the form
\begin{equation}
\label{NumberEight}
H_{\theta} = M + M(M^{-1} L)^{1-\theta}, \hspace{0.5in} \theta\in [0,1],
\end{equation}
are shown to be equivalent to their continuous counterparts on $\Lambda_{\theta}$, where $M$ and $L$ denote the mass and Laplacian matrices respectively assembled on the interface $\Gamma$. One particular example corresponds to $\theta=1/2$, which is shown in \cite{Arioli09} to adhere to the same coercivity and continuity bounds as the discrete Steklov-Poincar\'e operator, leading to mesh independent performance of GMRES. The norms presented in (\ref{NumberEight}) can be shown to be spectrally equivalent to
\[
\widetilde{H}_{\theta} = M(M^{-1} L)^{1-\theta}.
\]
Both of the above can be applied component-wise to a system, suggesting an appropriate form of $\widetilde{S}$ as 
\begin{equation}
\label{NumberNine}
\widehat{H}_{\theta} = \bigoplus_{1}^d \widetilde{H}_{\theta}.
\end{equation}
From the above, it is clear that fractional powers of matrices must be determined in order to apply the discrete norms. For relatively small problems, this can be achieved through direct methods such as a generalised eigenvalue decomposition. However, the complexity involved is $\mathcal{O}\left(n_{\Gamma}^3\right)$ suggesting instead the use of iterative approaches for larger problems. In \cite{Arioli09}, approximations through the use of truncated Lanczos and inverse Lanczos algorithms are described, and will also be employed within this work through the use of flexible GMRES \cite{Saad93} to account for the changing nature of the preconditioner.

\section{Results}

We present various results in this section to illustrate our approach in practice. It should be noted that while certain examples involve symmetric system matrices, our choice of non-symmetric preconditioner suggests GMRES as an appropriate Krylov solver.

\begin{figure}[!t]
\vspace{-1.3in}
\hspace{0.4in}
\psfrag{F}[c][r][0.8][0]{$\mathbf{g}$}
\includegraphics[width=0.9\textwidth]{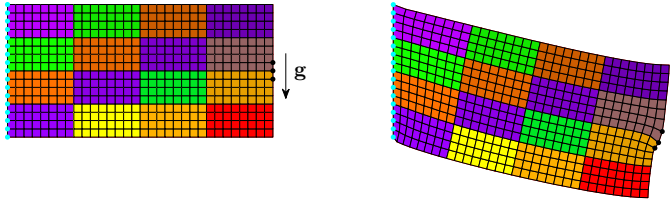}
\caption{Original decomposed undeformed layout (left) and deformed layout (right). Location of clamping and traction illustrated by cyan and black nodes, respectively.}
\label{DeflectionFigure}

\end{figure}

The test problem considered involves a cantilever beam over the 2D domain $\Omega = (0,2) \times (0,1)$, with downward force $\mathbf{f}=0.75$ and outward traction $\mathbf{g}=1$. The domain will be clamped on the right hand side through the application of homogeneous Dirichlet conditions on the relevant boundary. An illustration of the deflection as well as a pictorial example of a division of the domain (into $16$ subdomains) is provided in Figure \ref{DeflectionFigure}.

As discussed in the previous section, iterative approaches will be used for the application of $\widehat{H}_{\theta}$. For this work, it was found that the inverse Lanczos approach delivered the most promising results, largely due to the relatively small number of basis vectors required to apply the discrete norms.
Table \ref{TableOne} illustrates our results for differing mesh parameters $h$ and subdomains. The column labelled $\widetilde{S}=\hat{I}$ illustrates results for both test problems in the absence of interface preconditioning. By reading this column from top to bottom (for each problem), we observe a logarithmic dependence on the number of GMRES iterations for increasing mesh parameters. Reading this column from left to right also suggests a logarithmic dependence on the number of subdomains.

In comparison, the column labelled $\widetilde{S}=\widehat{H}_{\theta}$ of the table provides results with the interface preconditioner as discussed in (\ref{NumberNine}). Here, it can be seen that the number of iterations are independent of the chosen mesh parameter. Whilst there is a logarithmic dependence on the number of iterations for an increasing number of subdomains, a direct comparison with the column labelled $\widetilde{S}=\hat{I}$ suggests that our preconditioning strategy provides significant savings in the number of iterations required for convergence.

The final column labelled $\widetilde{S}=\widehat{H}_{\text{OPT}}$ illustrates results for selected values of theta based on testing. It was found that the recorded values were able to provide improved results over the other two columns, suggesting that different values of theta are able to provide a closer approximation to the decay of the associated Steklov-Poincar\'e operator.

\begin{table}[!h]
\begin{center}
{\centering{
{\resizebox{\textwidth}{!}{
\def\arraystretch{1.3}
\renewcommand{\tabcolsep}{0.5cm}
\begin{tabular}{ | r || c  c  c || c  c  c || c  c  c | }
\cline{2-10}
\multicolumn{1}{c|}{} & \multicolumn{3}{c||}{$\tilde{S} = I$} & \multicolumn{3}{c||}{$\tilde{S} = \widehat{H}_{\theta}$} & \multicolumn{3}{c|}{$\tilde{S} = \widehat{H}_{\text{OPT}}$} \\
\hline 
\, $\#$ Domains \, & 4 & 16 & 64 & 4 & 16 & 64 & 4 & 16 & 64 \\   
\hdashline
\, $\theta$ \, & - & - & - & 0.5 & 0.5 & 0.5 & 0.5 & 0.6 & 0.7 \\
\hline   
\hline
h = 1/32 \,\, & 28 & 47 & 68  & 12 & 18 & 27 & 12 & 17 & 22 \\
1/64 \,\, & 41 & 66 & 96 &  12 & 19 & 27 & 12 & 18 & 23  \\
1/128 \,\, & 59 & 93 & 137 &  12 & 19 & 27 & 12 & 19 & 24  \\ 
\hline  
\end{tabular} 
}}}}
\vspace{0.03in}
\caption{Results for the test problem. Tolerance of GMRES set at $10^{-6}$.}
\label{TableOne}
\end{center}
\end{table}
\vspace{-0.1in}

In order to observe the computational benefits of our method, we look to provide rough estimates in order to gauge how our derived approach will perform in a parallel environment. Due to the non-overlapping nature of our approach, all subdomain solves can be carried out in parallel. As mentioned previously, the main issue surrounds the solution to the resulting interface problem. Within each application of our preconditioner to this problem, we are required to invert the discrete interface Laplacian. This issue is present in the Lanczos process, and also in the subsequent generalised eigenvalue decomposition that follows. Due to the structure of this matrix, these inversions can lead to a computational bottleneck for an increasing number of subdomains, and so we would like to consider an iterative approach to alleviate this issue.

The structure of the involved matrix suggests conjugate gradient as a suitable alternative, coupled with an appropriate preconditioning strategy (PCG). In this work, we propose to precondition by using the relevant contributions of $L$ restricted to $\Gamma_i$, with the cross points removed to enable construction in parallel. The parallel CPU time taken for each GMRES iteration can then be realised by dividing the number of PCG iterations multiplied by the CPU time taken to apply the preconditioner by the total number of faces involved in the construction of $\Gamma$. By adding this contribution to the CPU time taken for one parallel subdomain solve, we calculate the total CPU time by multiplying the result to the total number of GMRES iterations required to achieve convergence. 

The results for the investigation are displayed in Table \ref{TableTwo} where CPU times (in seconds) are provided for differing mesh and subdomain sizes. A Linux machine with an Intel{\tiny{\textregistered}} Core\texttrademark \, i7 CPU 870 $@$ 2.93 GHz with 8 cores was used to obtain the data.  

\begin{table}[!h]
\begin{center}
{\centering{
{\resizebox{\textwidth}{!}{
\def\arraystretch{1.3}
\renewcommand{\tabcolsep}{0.8cm}
\begin{tabular}{ | r || c  c  c  c | }
\cline{2-5}
\multicolumn{1}{r|}{} & \multicolumn{4}{c|}{$\tilde{S} = \widehat{H}_{\text{OPT}}$} \\
\hline 
\, $\#$ Domains \, & 4 & 16 & 64 & 256 \\   
\hdashline
\, $\theta$ \, & 0.5 & 0.6 & 0.7 & 0.75 \\
\hline   
\hline
h = 1/16 & 0.0169 & 0.0168 & 0.0176  & 0.0254 \\
1/32 & 0.0635 & 0.0273 & 0.0199  & 0.0232   \\
1/64 & 0.4455 & 0.1238 & 0.0384  & 0.0238  \\ 
1/128 & 3.8716 & 1.0804 & 0.2529  & 0.0623   \\
1/256 & 50.2476 & 13.2858 & 3.3204  & 0.7295  \\ 
\hline  
\end{tabular} 
}}}}
\vspace{0.03in}
\caption{Total CPU times (seconds) anticipated through the use of parallel computing.}
\label{TableTwo}
\end{center}
\end{table}
\vspace{-0.1in}

From the table, it can be seen that for relatively coarse meshes, we do not see a significant enough decrease in the CPU time to warrant the use of parallelism. This behaviour can be attributed to the computational complexity of sparse matrix inversion ($\mathcal{O}(k^2 n)$, $k$ is the bandwidth) for relatively small values of $n$, and also the efficiency of the backslash command in MATLAB. However, notable savings in CPU time equating to roughly factor $4$ can be seen for finer meshes. These figures are encouraging, as they suggest that our approach is capable of significant speedup through the use of parallel architecture when compared directly to solving the problem globally on a single processor.

After collating the results in Table \ref{TableTwo}, a general increase was noted in the number of GMRES iterations when compared directly to the figures obtained in Table \ref{TableOne}. The reason for this can be attributed to the use of inner PCG iterations. In particular, a logarithmic dependence on the mesh parameter was observed for cases involving smaller numbers of subdomains. However, the deterioration can be seen as an acceptable compromise, as the results for larger meshes suggest the use of an increasing number of subdomains for improved performance. It should be noted that the results obtained above were done so with a relatively coarse tolerance for PCG of $10^{-3}$, as well as a reasonably modest number of PCG iterations (typically between $2$ and $12$) at each GMRES iteration for each of the test cases considered. 
 
It should not be expected that continual speedup can be gained through the use of an increasing number of subdomains, as certain factors such as inter-processor communication between each of the three steps will begin to play an important role. Therefore, in terms of a regular subdivision, this would suggest an optimal decomposition of the domain based on the mesh parameter, and also possibly other contributing factors relating to computer hardware.

\section{Topology Optimisation}

As an application of our findings, we will describe how our work can be incorporated into commonly used solvers from problems arising in topology optimisation. The problem we consider here is the so-called Variable Thickness Sheet problem \cite[pp.~54~--~57]{Bendsoe03}, which can be described mathematically using finite elements by the following nonlinear optimisation problem 
\vspace{-0.1in}
\begin{spacing}{1.1}
\begin{align}
\begin{cases}
\begin{aligned}
   &\min_{\mathbf{u},\boldsymbol{\rho}} && \mathbf{f}^{T} \mathbf{u}    \vspace{-0.3in} &&&& \\
   &\text{subject to:} && {K}(\boldsymbol{\rho} )\mathbf{u} =\mathbf{f} &&\hspace{0.4in} \left( K(\boldsymbol{\rho}) = \displaystyle\sum_{i=1}^m \rho_i K_i \right), \label{NumberTen} && \\
                        &&&\displaystyle\sum_{i \in D} \rho_{i} \leq V , &&&& \\
                        &&& 0 \leq \underline{\rho} \leq \rho_{i} \leq \overline{\rho} && \hspace{0.4in} \forall i \in D \defeq \left\{1,2,\dots,m\right\}. &&
\end{aligned}
\end{cases}
\end{align}
\end{spacing}
\vspace{-0.01in}
In the above, the $n-$vector of nodal displacement values $\mathbf{u} = \mathbf{u}(\boldsymbol{\rho})$ denotes the solution to the elasticity equations, with $\mathbf{f} \in \mathbb{R}^n$  representing the corresponding discretisation of the load linear form. The density $\boldsymbol{\rho} \in \mathbb{R}^m$ is subject to upper and lower bounds $\overline{\rho}$ and $\underline{\rho}$ respectively, with the volume of the body being denoted by $V$. Additionally, $K(\boldsymbol{\rho})$ represents the finite element stiffness matrix for the elasticity equations, with each $K_i$, $i=1, \dots, m$ denoting elemental stiffness matrices. 

By considering the method of Lagrange multipliers, minima to (\ref{NumberTen}) are obtained through a nonlinear system of equations. The nonlinearities can be dealt with using a number of commonly used approaches. For instance, one could consider the use of interior point methods \cite{Hoppe04}. The fairly standard solution technique used by the community involves the consideration of fixed point type update schemes for an initial guess for the density in the following way
\vspace{0.1in}
\begin{enumerate}
\item Finite Element Analysis (FEA) – solve equations of linear elasticity.
\item Density update – e.g.: OC (see \cite{Bendsoe03}), MMA (see \cite{Svanberg02}).
\item Check for convergence. If not satisfied, rerun 1 and 2 using updated density.
\end{enumerate}
\vspace{0.1in}
It can be expected that a reasonably large number of fixed point iterations are required to obtain a suitable final design. The bulk of computational effort will be concentrated on the Finite Element Analysis step, namely the repeated process of obtaining updated displacement variables through the use of the equations of linear elasticity \cite{Borrvall01}. Therefore, we propose to apply our preconditioning strategy as discussed in Section \ref{SectionThree} to this problem coupled with the fairly straightforward Optimality Criteria (OC) method for the density update. No attempt will be made here to carry out Step $2$ above in parallel; however \cite{Borrvall01} describe an appropriate implementation using the Method of Moving Asymptotes (MMA).

In Table \ref{TableThree}, results are provided illustrating the performance of our approach for the cantilever beam problem. The results were obtained using an adaptive tolerance for GMRES based on successive compliance values. The total number of fixed point iterations are given, along with the average number of GMRES iterations per fixed point step (bracketed). Whilst the number of fixed point iterations appears to increase for finer meshes, the average number of GMRES iterations remains roughly constant. Whilst we still see a logarithmic dependence on the average number of GMRES iterations for an increasing number of subdomains, the fixed point iterations remain roughly constant. 

Future work involves validation of our approach on a parallel machine, as well as consideration of further problems (possibly to include 3D domains) and alternative solution methods to try to solve topology optimisation problems completely in parallel. We expect our approach to adapt well in parallel, with potential speedup for 3D problems of factor $8$ anticipated.
\begin{table}[!t]
{\centering{
{\resizebox{\textwidth}{!}{
\def\arraystretch{1.3}
\renewcommand{\tabcolsep}{0.8cm}
\begin{tabular}{ | r || c  c  c  c | }
\hline 
\, $\#$ Domains \, & 4 & 16 & 64 & 256 \\   
\hdashline
\, $\theta$ \, & 0.5 & 0.6 & 0.7 & 0.75 \\
\hline   
\hline
h = 1/16 \,\, & 10 (10) & 10 (18) & 10 (33)  & 10 (56) \\
1/32 \,\, & 17 (11) & 17 (18) & 17 (34)  & 19 (54)   \\
1/64 \,\, & 23 (11) & 23 (18) & 24 (32)  & 27 (54)  \\ 
1/128 \,\, & 29 (12) & 30 (17) & 32 (31)  & 32 (52)   \\
1/256 \,\, & 33 (13) & 36 (17) & 37 (31)  & 41 (49)  \\ 
\hline
\end{tabular} 
}}}}
\vspace{0.03in}
\centering
\caption{Results for the cantilever beam problem solved using our preconditioning strategy for the FEA coupled with the OC method for the density update.} 
\label{TableThree}
\vspace{-0.2in}
\end{table}

\vspace{-0.05in}

\bibliographystyle{siam}
\bibliography{myrefs1}

\end{document}